\documentclass[12pt, reqno]{amsart}
\usepackage{ amsmath,amsthm, amscd, amsfonts, amssymb, graphicx, color}
\usepackage[bookmarksnumbered, colorlinks, plainpages]{hyperref}
\newtheorem{thm}{Theorem}[section]
\newtheorem{cor}[thm]{Corollary}
\newtheorem{lem}[thm]{Lemma}

\newtheorem{exam}[thm]{Example}
\numberwithin{equation}{section}

\begin{document}

\title{Study on Hirano Invertibility of Two New Types of Perturbed Operator Matrices }

\author[Gou Haibo]{Gou Haibo}
\address{Gou Haibo\\ School of Mathematics \\ Hangzhou Normal University \\ China}
\email{<haibo\_gou@163.com>}

\author[Chen Huanyin]{Chen Huanyin}
\address{Chen Huanyin \\ School of Mathematics \\ Hangzhou Normal University \\ China}
\email{<huanyinchen@aliyun.com>}
\thanks{}


\begin{abstract}
We investigate the Hirano invertibility of block-operator matrices in Banach algebras, and obtain the Hirano inverse of matrix $\begin{bmatrix}
A&B\\
C&D
\end{bmatrix}$ under two types of new perturbation conditions. Furthermore, we provide a new operator matrix decomposed into the sum of tripotent and nilpotent elements on Banach spaces.
\end{abstract}

\maketitle

\section{Introduction}

Let $\mathcal L(X)$ be a Banach algebra of all bounded linear operators over a Banach space $X$. An operator $A\in \mathcal L(X)$ has Drazin inverse if there exists $Z\in \mathcal L(X)$ such that $$AZ=ZA, Z=ZAZ, A-A^2Z\in \mathcal L(X)$$ is nilpotent. Such $Z,$ if it exists, is unique and will be denoted by $A^D$ (see [1, Theorem 2.4]).

Recall that an operator  $A\in \mathcal L(X)$ has strongly Drazin inverse if there exists $Z\in \mathcal L(X)$ such that $$AZ=ZA, Z=ZAZ, A-AZ\in \mathcal L(X)$$ is nilpotent. Such $Z,$ if it exists, is unique and will be denoted by $A^{sD}$ (see [2, Theorem 2.1]).

An operator $A\in \mathcal L(X)$ has Hirano inverse if there exists $Z\in \mathcal L(X)$ such that $$AZ=ZA, Z=ZAZ, A^2-AZ\in \mathcal L(X)$$ is nilpotent. Such $Z,$ if it exists, is unique and will be denoted by $A^H$ (see [3, Corollary 2.2]).

Evidently, $A$ has strongly Drazin inverse if and only if it is the sum of an idempotent and a nilpotent operators that commute. $A$ has Hirano inverse if and only if it is the sum of a tripotent and a nilpotent operators that commute (see [3, Theorem 3.1]).

The following theorem gives the characteristics of Hirano invertibility:

\begin{thm}
 If $A\in \mathcal L(X)$, the following conditions are equivalent:
 \begin{itemize}
 \item $A\in \mathcal L(X)$ has Hirano inverse.
 \item $A-A^3\in N(\mathcal L(X))$.
 \item $A$ is the sum of commutative tripotent and nilpotent operators.
 \item $A^2$ has strongly Drazin inverse.
 \end{itemize}
\end{thm}

Note that all complex matrices of order $n$ are bounded linear operators on Banach spaces. In particular, if $A$ is a complex matrix of order $n$, then $A$ has a Hirano inverse if and only if its eigenvalues are $-1,0$ or $1$. In this paper, we study the Hirano inverse of block-operator matrix $\begin{bmatrix}
A&B\\
C&D
\end{bmatrix}$ under perturbation conditions, and provide a new operator matrix that can be decomposed into the sum of tripotent and nilpotent operators.

In section 2, we study the Hirano inverse of operator matrix $\begin{bmatrix}
A&B\\
C&D
\end{bmatrix}$ with strongly Drazin invertible subblocks. In section 3, we study the Hirano inverse of operator matrices with Hirano invertible subblocks.  A new perturbation condition is obtained for two types of block-operator matrices to be expressed as the sum of commutative tripotent and nilpotent operators.

In this paper, $\mathcal L(X)$ represents a Banach algebra composed of all bounded linear operators on a Banach space $X$. $N(\mathcal L(X))$ represents the nilpotent element in $\mathcal L(X)$, $M_n(\mathcal L(X))$ represents matrices with order $n$ in $\mathcal L(X)$. Set $A^e=AA^D, A^{\pi}=I-A^e$. We call the smallest positive integer $k$ that satisfies $(A-A^2A^D)^k=0$ an index of $A$ and record it as $ind(A)=k$. 

Note $l^2(\mathbb N)=(x=\{x_k\}: \sum\limits_{k=1}\limits^{\infty}{x_k^2\textless \infty})$.

\section{Operator matrices with strongly Drazin invertible subblocks}

The aim of this section is to investigate the Hirano invertibility of operator matrix $\begin{bmatrix}
A&B\\
C&D
\end{bmatrix}$ with strongly Drazin invertible subblocks. Firstly, some lemmas for the strongly Drazin inverse and Hirano inverse needed in the following article are given.

\begin{lem}
Let $P, Q\in \mathcal L(X)$ and $PQ=0$. If $P$ and $Q$ have strongly Drazin inverses, then $P+Q$ has strongly Drazin inverse.
\end{lem}

\begin{proof} 
Since $P$ and $Q$ have strongly Drazin inverses, we have $P-P^2, Q-Q^2\in N(\mathcal L(X)).$

Clearly,
$$P+Q=\begin{bmatrix}
1&Q\\
\end{bmatrix}\begin{bmatrix}
P\\
1
\end{bmatrix}, \begin{bmatrix}
P\\
1
\end{bmatrix}\begin{bmatrix}
1&Q\\
\end{bmatrix}=\begin{bmatrix}
P&0\\
1&Q
\end{bmatrix}.$$
$$
\begin{bmatrix}
P&0\\
1&Q
\end{bmatrix}-\begin{bmatrix}
P&0\\
1&Q
\end{bmatrix}^2=\begin{bmatrix}
P-P^2&0\\
1-(P+Q)&Q-Q^2
\end{bmatrix}\in N(\mathcal L(X)).
$$

By [2, Theorem 2.1], $\begin{bmatrix}
P&0\\
1&Q
\end{bmatrix}$is strongly Drazin invertible. According to Cline's formula (see[4, Lemma 2.1]), $P+Q$ is also strongly Drazin invertible.
\end{proof}

\begin{lem}
Let $P, Q\in \mathcal L(X)$ and $PQP=0, PQ^2=0.$ If $P$ and $Q$ have Hirano inverses, then $P+Q$ has Hirano inverse.
\end{lem}

\begin{proof}
We compute that
$$
P+Q=\begin{bmatrix}
1&Q\\
\end{bmatrix}\begin{bmatrix}
P\\
1
\end{bmatrix}, \begin{bmatrix}
P\\
1
\end{bmatrix}\begin{bmatrix}
1&Q\\
\end{bmatrix}=\begin{bmatrix}
P&PQ\\
1&Q
\end{bmatrix}.$$
$$
M:=\begin{bmatrix}
P&PQ\\
1&Q
\end{bmatrix}^2=\begin{bmatrix}
P^2+PQ&P^2Q\\
P+Q&PQ+Q^2
\end{bmatrix}.
$$

Then
$$M=\begin{bmatrix}
PQ&P^2Q\\
0&PQ
\end{bmatrix}+\begin{bmatrix}
P^2&0\\
P+Q&Q^2
\end{bmatrix}=C+D.
$$

We check that $CD=0, C\in N(\mathcal L(X)),$ 

$$
D-D^2=\begin{bmatrix}
P(P-P^3)&0\\
(P+Q)(1-P^2)+Q^2(P+Q)&Q(Q-Q^3)
\end{bmatrix}\in N(\mathcal L(X)).
$$

Thus, D is strongly Drazin invertible. Therefore, by Lemma 2.1, $M$ is  strongly Drazin invertible. According to [4, Theorem 2.4], an operator $X$ has Hirano inverse if and only if $X^2$ has strongly Drazin inverse in a Banach algebra.

Consequently, $P+Q$ is Hirano invertible.
\end{proof}

\begin{lem}
Let $A,B\in \mathcal L(X)$ have strongly Drazin inverses and $A^DBA^D=0$, then 
$\begin{bmatrix}
AA^e&B\\
A^e&0
\end{bmatrix}$ has strongly Drazin inverse.
\end{lem}

\begin{proof}
We have matrix decomposition:
$$\begin{bmatrix}
AA^e&B\\
A^e&0
\end{bmatrix}=\begin{bmatrix}
AA^e&I\\
A^e&0
\end{bmatrix}\begin{bmatrix}
I&0\\
0&B
\end{bmatrix},
$$

$$\begin{bmatrix}
I&0\\
0&B
\end{bmatrix}\begin{bmatrix}
AA^e&I\\
A^e&0
\end{bmatrix}=\begin{bmatrix}
AA^e&I\\
BA^e&0
\end{bmatrix}.
$$
Obviously, $\begin{bmatrix}
AA^e&I\\
BA^e&0
\end{bmatrix}=\begin{bmatrix}
AA^e&0\\
0&0
\end{bmatrix}+\begin{bmatrix}
0&I\\
BA^e&0
\end{bmatrix}:=C+D$.

Besides,
$$C-C^2=\begin{bmatrix}
(A-A^2)A^e&0\\
0&0
\end{bmatrix}, D-D^2=\begin{bmatrix}
-BA^e&I\\
BA^e&-BA^e
\end{bmatrix}.
$$

Because $A$ has strongly Drazin inverse, then $A-A^2\in N( \mathcal L(X)),$ thus $C-C^2\in N( \mathcal L(X)).$ Clearly, $(D-D^2)^4=0.$ Accordingly, $C,D$ have strongly Drazin inverses and satisfy $CDC^2=0, CDCD=0, CD^2=0.$ By[2, Theorem 2.7], $\begin{bmatrix}
AA^e&I\\
BA^e&0
\end{bmatrix}$has strongly Drazin inverse. By use of Cline's formula,  $\begin{bmatrix}
AA^e&B\\
A^e&0
\end{bmatrix}$ has strongly Drazin inverse.
\end{proof}

\begin{lem}
Let $M=\begin{bmatrix}
A&B\\
I&0
\end{bmatrix}$ and $A,B$ have strongly Drazin inverses. If $A^DBA^D=0, BA^{\pi}B=0$ and $BAA^{\pi}=0,$ then $M$ has Hirano inverse.
\end{lem}

\begin{proof}
Obviously, we have
$$M=\begin{bmatrix}
AA^e&B\\
A^e&0
\end{bmatrix}+\begin{bmatrix}
AA^{\pi}&0\\
A^{\pi}&0
\end{bmatrix}=P+Q.$$

According to Lemma 2.3, $P$ is strongly Drazin invertible and then $P$ is Hirano invertible. Set $k=ind(A).$ Then $(A-A^2A^D)^k=0,$ and so $A^kA^{\pi}=0.$

Clearly, we have
$$
Q^i=\begin{bmatrix}
A^iA^{\pi}&0\\
A^{i-1}A^{\pi}&0
\end{bmatrix}\ 
i=1,2,3...
$$

Then $(Q-Q^3)^k=0,$ thus $Q$ is Hirano invertible. 

Clearly, $PQP=0, PQ^2=0.$ According to Lemma 2.3, $M$ has Hirano inverse.
\end{proof}

\begin{lem}
Let $M=\begin{bmatrix}
A&B\\
C&0
\end{bmatrix}$ and $A, BC$ have strongly Drazin inverses. If $A^DBCA^D=0, BCA^{\pi}BC=0$ and $BCA^{\pi}A=0,$ then $M$ has Hirano inverse.
\end{lem}

\begin{proof}
Observing that 
$$
\begin{bmatrix}
A&B\\
C&0
\end{bmatrix}=\begin{bmatrix}
I&0\\
0&C
\end{bmatrix}\begin{bmatrix}
A&B\\
I&0
\end{bmatrix}\ and
$$

$$
\begin{bmatrix}
A&B\\
I&0
\end{bmatrix}\begin{bmatrix}
I&0\\
0&C
\end{bmatrix}=\begin{bmatrix}
A&BC\\
I&0
\end{bmatrix}.
$$
Applying Lemma 2.4 and Cline's formula, we obtain the result.
\end{proof}

\begin{lem}
Let $A, D, BC$ have strongly Drazin inverses. If $ABC=0, BCA^{\pi}=0, BDC=0$ and $BD^2=0,$ then $\begin{bmatrix}
A&B\\
C&D
\end{bmatrix}$ has Hirano inverse.
\end{lem}

\begin{proof}
Clearly, we have
$$\begin{bmatrix}
A&B\\
C&D
\end{bmatrix}=\begin{bmatrix}
A&B\\
C&0
\end{bmatrix}+\begin{bmatrix}
0&0\\
0&D
\end{bmatrix}=P+Q.
$$

Apply Lemma 2.5, $P$ has Hirano inverse. Since $D$ has strongly Drazin inverse, then $Q$ has Hirano inverse. As $PQP=0$ and $PQ^2=0,$ $P+Q$ has Hirano inverse.
\end{proof}

\begin{thm}
Let $M=\begin{bmatrix}
A&B\\
C&D
\end{bmatrix}, A$ and $D$ have strongly Drazin inverses. If $BDD^D=0, D^{\pi}CB=0, D^{\pi}CA=0,$ then $M$ has Hirano inverse.
\end{thm}

\begin{proof}
We verify that
$$M=\begin{bmatrix}
0&BDD^D\\
D^{\pi}C&0
\end{bmatrix}+\begin{bmatrix}
A&BD^{\pi}\\
DD^DC&D
\end{bmatrix}=P+Q.
$$

Since
$P^2=0, P$ has Hirano inverse. 

Apply Lemma 2.6 to $Q$, it follows from $$ABD^{\pi}(DD^DC)=0, BD^{\pi}D^2=0\ and\ BD^{\pi}D(DD^DC)=0,$$ $Q$ has Hirano inverse.
Since 
$$
PQP=0, PQ^2=0,
$$
$M=P+Q$ has Hirano inverse as asserted.
\end{proof}

\begin{cor}
Let $M=\begin{bmatrix}
A&B\\
C&D
\end{bmatrix}, A$ and $D$ have strongly Drazin inverses. If $BD=0, D^{\pi}C=0,$ then $M$ has Hirano inverse.
\end{cor}

\begin{proof}
By virtue of Theorem 2.7, the proof is true.
\end{proof}

\begin{cor}
Let $D$ has strongly Drazin inverse. If $CB=0, BD=0$ and $CA=0,$ then $\begin{bmatrix}
A&B\\
C&D
\end{bmatrix}$ has Hirano inverse.
\end{cor}

\begin{proof}
By Theorem 2.7, the result is true.
\end{proof}

\begin{thm}
Let $M=\begin{bmatrix}
A&B\\
C&D
\end{bmatrix}, A$ and $D$ have strongly Drazin inverses. If $A^DBCA^D=0, A^DBD=0, A^{\pi}BC=0$ and $A^{\pi}BD=0$ then $M$ has Hirano inverse.
\end{thm}

\begin{proof}
Observing that 
$$M=\begin{bmatrix}
A&AA^DB\\
C&D
\end{bmatrix}+\begin{bmatrix}
0&A^{\pi}B\\
0&0
\end{bmatrix}=P+Q.
$$

We have
$$
PQ^2=\begin{bmatrix}
0&AA^{\pi}B\\
0&CA^{\pi}B
\end{bmatrix}
\begin{bmatrix}
0&A^{\pi}B\\
0&0
\end{bmatrix}=0,
$$

$$
PQP=\begin{bmatrix}
0&AA^{\pi}B\\
0&CA^{\pi}B
\end{bmatrix}
\begin{bmatrix}
A&AA^DB\\
C&D
\end{bmatrix}=\begin{bmatrix}
AA^{\pi}BC&AA^{\pi}BD\\
CA^{\pi}BC&CA^{\pi}BD
\end{bmatrix}=0.
$$

Write
$$
P=\begin{bmatrix}
A^2A^D&AA^DB\\
CAA^D&D
\end{bmatrix}+\begin{bmatrix}
AA^{\pi}&0\\
CA^{\pi}&0
\end{bmatrix}=P_1+P_2.
$$

By virtue of  Lemma 2.6, we prove that $P_1$ has Hirano inverse.

Computing that 
$$
\begin{array}{r l l}
(P_2-P_2^3)^3&=&\begin{bmatrix}
(A-A^3)^3A^{\pi}&0\\
(C-CA^2)(A-A^3)^2A^{\pi}&0
\end{bmatrix}\\
&=&\begin{bmatrix}
(A-A^3)(A+A^2-A^3-A^4)A^{\pi}(A-A^2)&0\\
(C-CA^2)A^{\pi}(A+A^2-A^3-A^4)A^{\pi}(A-A^2)&0
\end{bmatrix}\\
&\in& N(\mathcal L(X)).
\end{array}
$$

Therefore, $P_2$ has Hirano inverse. We directly verify that
$$
P_2P_1=0,
$$

By Lemma 2.1, $P=P_1+P_2$ has Hirano inverse, thus $M$ has Hirano inverse.
\end{proof}

\begin{cor}
Let $M=\begin{bmatrix}
A&B\\
C&D
\end{bmatrix}$ and $A, D$ have strongly Drazin inverses. If $D^DCBD^D=0, D^DCA=0, D^{\pi}CB=0$ and $D^{\pi}CA=0,$ then $M$ has Hirano inverse.
\end{cor}

\begin{proof}
Let $P=\begin{bmatrix}
0&I\\
I&0
\end{bmatrix}$, then $M=P\begin{bmatrix}
D&C\\
B&A
\end{bmatrix}P.$ $M$ has Hirano inverse if and only if $\begin{bmatrix}
D&C\\
B&A
\end{bmatrix}$ has Hirano inverse (see[14, Theorem 3.6]).

Write
$$\begin{bmatrix}
D&C\\
B&A
\end{bmatrix}=\begin{bmatrix}
D&DD^DC\\
B&A
\end{bmatrix}+\begin{bmatrix}
0&D^{\pi}C\\
0&0
\end{bmatrix}=P+Q.
$$

Computing that
$$
PQP=\begin{bmatrix}
0&DD^{\pi}C\\
0&BD^{\pi}C
\end{bmatrix}
\begin{bmatrix}
D&DD^DC\\
B&A
\end{bmatrix}=\begin{bmatrix}
DD^{\pi}CB&DD^{\pi}CA\\
BD^{\pi}CB&BD^{\pi}CA
\end{bmatrix}=0,
$$

$$
PQ^2=\begin{bmatrix}
0&DD^{\pi}C\\
0&BD^{\pi}C
\end{bmatrix}\begin{bmatrix}
0&D^{\pi}C\\
0&0
\end{bmatrix}=0.
$$

Observing that 
$$
P=\begin{bmatrix}
D^2D^D&DD^DC\\
BDD^D&A
\end{bmatrix}+\begin{bmatrix}
DD^{\pi}&0\\
BD^{\pi}&0
\end{bmatrix}=
P_1+P_2.
$$

By virtue of Lemma 2.6, we prove that $P_1$ has Hirano inverse. $P_2-P_2^3 \in N(\mathcal L(X)), P_2$ has Hirano inverse. We directly verify that
$$
P_2P_1=0,
$$

By Lemma 2.2, $P=P_1+P_2$ has Hirano inverse, thus $M$ has Hirano inverse.
\end{proof}

Assuming that $I$ is an identical linear operator, we give a numerical example of the above theorem.

\begin{exam}
Let $A, B, C, D$ be bounded linear operators acting on $l^2(\mathbb N).$ $M=\begin{bmatrix}
A&B\\
C&D
\end{bmatrix}$, defined as follows respectively:
$$A(x_1, x_2, x_3, ...)=(x_1+x_3, x_2+x_3, x_3,...),$$
$$B(x_1, x_2, x_3, ...)=(x_3, x_3, -x_3,...),$$
$$C(x_1, x_2, x_3, ...)=(x_1, x_1+x_2, 0,...),$$
$$D(x_1, x_2, x_3, ...)=(x_1+x_3,x_2+x_3,0,...).$$

Because $(A-A^2)(x_1, x_2, x_3, ...)=(-x_3, -x_3, 0,...)\in N(\mathcal L(X)),$ thus $A$ has strongly Drazin inverse. Obviously, $D$ is an idempotent operator, and by the same token, $D$ has strongly Drazin inverse. Then $D^D=D, D^{\pi}(x_1, x_2, x_3, ...)=(-x_3, -x_3, x_3,...).$ Furthermore, $BDD^D=0, D^{\pi}CB=0, D^{\pi}CA=0.$ By Theorem 2.7, $M$ has Hirano inverse. In fact, the eigenvalues of $M$ are 0,1 (quintuple).
\end{exam}

\section{Operator matrices with Hirano invertible subblocks}

The purpose of this section is to investigate the Hirano invertibility of block-operator matrix $\begin{bmatrix}
A&B\\
C&D
\end{bmatrix}$. We first prove the new additive properties of the perturbed operator Hirano inverse.

\begin{lem}
Let $A \in N(\mathcal L(X))$ and $B\in \mathcal L(X)$ has Hirano inverse. If $AB^H=0, B^{\pi}AB=0,$ then $A+B$ has Hirano inverse.
\end{lem}

\begin{proof} 
Let $P=BB^H.$ Then $B$ can be represented as
$$
\begin{array}{r l l}
B&=&\begin{bmatrix}
PBP&PB(1-P)\\
(1-P)BP&(1-P)B(1-P)
\end{bmatrix}_P\\
&=&\begin{bmatrix}
B^2B^H&0\\
0&BB^{\pi}
\end{bmatrix}_P.
\end{array}
$$ 

Set 
$A=\begin{bmatrix}
A_1&A_3\\
A_4&A_2
\end{bmatrix}_P.$ From $AB^H=0$ it follows that $A=\begin{bmatrix}
0&A_3\\
0&A_2
\end{bmatrix}_P,$ where $A_2$ is nilpotent.

Then 
$$
\begin{array}{r l l}
A+B&=&\begin{bmatrix}
B^2B^H&A_3\\
0&BB^{\pi}+A_2
\end{bmatrix}_P\\
&=&\begin{bmatrix}
B^2B^H&A_3\\
0&BB^{\pi}+B^{\pi}AB^{\pi}
\end{bmatrix}_P
\end{array}
$$
 
 Since $A=AB^{\pi}$ has Hirano inverse, it follows that $(AB^{\pi})B^{\pi}$ has Hirano inverse. By using Cline's formula, 
$B^{\pi}(AB^{\pi})$ has Hirano inverse. 

That is, $B^{\pi}AB^{\pi}$ has Hirano inverse. Clearly, $BB^{\pi}$ is nilpotent, and so it has Hrano inverse.

Since $(B^{\pi}AB^{\pi})(BB^{\pi})=B^{\pi}AB=0$, 
we prove that  $B^{\pi}AB^{\pi}+BB^{\pi}=BB^{\pi}+A_2$ has Hirano inverse, as required.

Additively, $B^2B^H$ has Hirano inverse, it follows that $A+B$ has Hirano inverse.

\end{proof}

\begin{lem}
Let $M=\begin{bmatrix}
A&0\\
C&D
\end{bmatrix}\in M_2( \mathcal L(X)).$ If $A, D$ have Hirano inverses, then $M\in M_2( \mathcal L(X))^H$.
\end{lem}

\begin{proof} 
By use of [2, Theorem 3.1], the result is obtained.
\end{proof}

\begin{lem}
Let $A, B\in \mathcal L(X)$ have Hirano inverses. If $A^HB=0, AB^H=0$ and $B^{\pi}ABA^{\pi}=0$, then $A+B\in \mathcal L(X)$ has Hirano inverse.
\end{lem}

\begin{proof}
Let $P=AA^H.$ Then $A$ can be represented as
$$
\begin{array}{r l l}
A&=&\begin{bmatrix}
PAP&PA(1-P)\\
(1-P)AP&(1-P)A(1-P)
\end{bmatrix}_P\\
&=&\begin{bmatrix}
A^2A^H&0\\
0&AA^{\pi}
\end{bmatrix}_P.
\end{array}
$$ 

Set 
$B=\begin{bmatrix}
B_1&B_2\\
B_3&B_4
\end{bmatrix}_P.$ From $A^HB=0$ it follows that 
$$
B=\begin{bmatrix}
0&0\\
B_3&B_4
\end{bmatrix}_P=\begin{bmatrix}
0&0\\
BAA^H&BA^{\pi}
\end{bmatrix}_P
$$

Then 
$$
A+B=\begin{bmatrix}
A^2A^H&0\\
B_3&AA^{\pi}+B_4
\end{bmatrix}_P
$$

From the assumptions $AB^H=0, B^{\pi}ABA^{\pi}=0,$ we have
$$
\begin{array}{r l l}
AB^H&=&\begin{bmatrix}
A^2A^H&0\\
0&AA^{\pi}
\end{bmatrix}_P\begin{bmatrix}
0&0\\
(B_4^H)^2B_3&B_4^H
\end{bmatrix}_P\\
&=&\begin{bmatrix}
0&0\\
(AA^{\pi})(B_4^H)^2B_3&(AA^{\pi})B_4^H
\end{bmatrix}_P\\
&=&0,
\end{array}
$$

$$
\begin{array}{r l l}
 B^{\pi}ABA^{\pi}&=&\begin{bmatrix}
1&0\\
-B_4(B_4^H)^2B_3&B_4^{\pi}
\end{bmatrix}_P\begin{bmatrix}
0&0\\
AA^{\pi}b_3&AA^{\pi}B_4
\end{bmatrix}_P\begin{bmatrix}
A^{\pi}&0\\
0&1
\end{bmatrix}_P\\
&=&\begin{bmatrix}
0&0\\
0&(B_4)^{\pi}(AA^{\pi})B_4
\end{bmatrix}_P\\
 &=&0.
 \end{array}
$$

Thus, $(AA^{\pi})(B_4)^H=0$ and $(B_4)^{\pi}(AA^{\pi})B_4=0.$

Using Lemma 3.1, we know $AA^{\pi}+B_4$ has Hirano inverse. Therefore, $A+B$ has Hirano inverse obtained.
\end{proof}

\begin{thm}
Let $M=\begin{bmatrix}
A&B\
C&D
\end{bmatrix}\in M_2( \mathcal L(X)), A, D\in \mathcal L(X)$ have Hirano inverses. If $BD^H=0, A^{\pi}BC=0, A^HBC=0, (DD^{\pi}-CA^HB)C=0,$ then $M\in M_2( \mathcal L(X)^H$.
\end{thm}

\begin{proof}
Write
$$
M=\begin{bmatrix}
0&B\\
0&DD^{\pi}
\end{bmatrix}+\begin{bmatrix}
A&0\\
C&D^2D^H
\end{bmatrix}.
$$

Clearly, 
$$
P=\begin{bmatrix}
0&B\\
0&DD^{\pi}
\end{bmatrix},
Q=\begin{bmatrix}
A&0\\
C&D^2D^H
\end{bmatrix},
$$
where
$$
P^2=\begin{bmatrix}
0&BDD^{\pi}\\
0&D^2D^{\pi}
\end{bmatrix}\in N(M_2(\mathcal L(X))).
$$

Thus, $P\in N(M_2(\mathcal L(X)).$

Obviously,
$$
Q-Q^3=
\begin{bmatrix}
A-A^3&0\\
C-T&D(D-D^3)D^H
\end{bmatrix}\in N(M_2( \mathcal L(X))),
$$
where $T=CA^2+D^2D^HCA+D^3D^HC$.

Therefore, $Q\in M_2( \mathcal L(X))^H$. In [2, Corollary 3.2], $(D^H)^H=D^2D^H.$

Observe that 
$$
PQ^H=
\begin{bmatrix}
0&B\\
0&DD^{\pi}
\end{bmatrix}\begin{bmatrix}
A^H&0\\
G&D^H
\end{bmatrix}=\begin{bmatrix}
BG&BD^H\\
DD^{\pi}G&0
\end{bmatrix}=0,
$$

$$
\begin{array}{r l l}
Q^{\pi}PQ&=&
\begin{bmatrix}
A^{\pi}&0\\
M&D^{\pi}
\end{bmatrix}\begin{bmatrix}
0&B\\
0&DD^{\pi}
\end{bmatrix}\begin{bmatrix}
A&0\\
C&D^2D^H
\end{bmatrix}\\
&=&\begin{bmatrix}
A^{\pi}BC&A^{\pi}BD^2D^H\\
MB+DD^{\pi}C&(MB+DD^{\pi})D^2D^H
\end{bmatrix}\\
&=&0,
\end{array}
$$
where
$$G=\sum\limits_{i=0}\limits^{r-1}{(D^H)^{i+2}CA^iA^{\pi}+\sum\limits_{i=0}\limits^{s-1}{D^{\pi}(D^2D^H)^iC(A^H)^{i+2}}-D^HCA^H},$$
$$ind(A)=r, ind(D)=s,$$
$$M=-(CA^H+D^2D^HG).$$

According to Lemma 3.1, $M=P+Q$ has Hirano inverse.
\end{proof}

\begin{cor}
Let $M=\begin{bmatrix}
A&B\\
C&D
\end{bmatrix}, A, D\in \mathcal L(X)$ have Hirano inverses. If $BD^H=0, BC=0,$ and $DD^{\pi}C=0,$ then $M\in M_2( \mathcal L(X))^H$.
\end{cor}

\begin{proof}
Compared with the previous conditions in Theorem 3.4, it is easier to prove $M\in M_2( \mathcal L(X))^H$.
\end{proof}

\begin{exam}
Let $M=\begin{bmatrix}
A&B\\
C&D
\end{bmatrix}$, where $A=\begin{bmatrix}
1&0\\
2&1
\end{bmatrix},D=\begin{bmatrix}
1&0\\
1&0
\end{bmatrix},\\ C=\begin{bmatrix}
0&1\\
0&1
\end{bmatrix}$ and $B=\begin{bmatrix}
1&-1\\
-1&1
\end{bmatrix}.$ Observing that
$$
A-A^3=\begin{bmatrix}
0&0\\
-4&0
\end{bmatrix}\in N(\mathbb C^{2\times 2}), D^2=D.
$$

Thus $A,D \in M_2(\mathbb C)^H.$
 $$
 BD^H=\begin{bmatrix}
1&-1\\
-1&1
\end{bmatrix}\begin{bmatrix}
1&0\\
1&0
\end{bmatrix}=0, BC=\begin{bmatrix}
1&-1\\
-1&1
\end{bmatrix}\begin{bmatrix}
0&1\\
0&1
\end{bmatrix}=0,
$$
$$
DD^{\pi}C=\begin{bmatrix}
1&0\\
1&0
\end{bmatrix}\begin{bmatrix}
0&0\\
-1&1
\end{bmatrix}\begin{bmatrix}
0&1\\
0&1
\end{bmatrix}=0.
$$ 
By Corollary 3.5, $M\in M_4( \mathbb C)^H$.
\end{exam}

\begin{thm}
Let $A, D\in \mathcal L(X)$ have Hirano inverses. If $AB=0, BD^H=0$ and $D^{\pi}CB=0$, then $M=\begin{bmatrix}
A&B\\
C&D
\end{bmatrix}\in M_2(\mathcal L(X))^H.$
\end{thm}

\begin{proof}
Let $M=P+Q,$ where 
$P=\begin{bmatrix}
A&0\\
C&0
\end{bmatrix}, Q=\begin{bmatrix}
0&B\\
0&D
\end{bmatrix}.$

According to [4, Proposition 4.1 and 4.2], 
$$
P^H=\begin{bmatrix}
A^H&0\\
C(A^H)^2&0
\end{bmatrix}, Q^H=\begin{bmatrix}
0&0\\
0&D^H
\end{bmatrix}.
$$

Notice that $A^HB=(A^H)^2AB=0,$ we have
$$
\begin{array}{r l l}
P^HQ=\begin{bmatrix}
A^H&0\\
C(A^H)^2&0
\end{bmatrix}\begin{bmatrix}
0&B\\
0&D
\end{bmatrix}&=&\begin{bmatrix}
0&A^HB\\
0&C(A^H)^2B
\end{bmatrix}\\
&=&\begin{bmatrix}
0&(a^H)^2AB\\
0&C(A^H)^3AB
\end{bmatrix}\\
&=&0,
\end{array}
$$

$$
PQ^H=\begin{bmatrix}
A&0\\
C&0
\end{bmatrix}\begin{bmatrix}
0&0\\
0&D^H
\end{bmatrix}=0,
$$

$$
Q^{\pi}PQ=\begin{bmatrix}
1&0\\
0&D^{\pi}
\end{bmatrix}\begin{bmatrix}
0&0\\
0&CB
\end{bmatrix}=0.
$$

By use of Lemma 3.3, $M=\begin{bmatrix}
A&B\\
C&D
\end{bmatrix}\in M_2(\mathcal L(X))^H.$
\end{proof}

\begin{cor}
Let $A,D\in \mathcal L(X)$ have Hirano inverses. If $CB=0, CA=0$ and $A^HBC=0$, then $M=\begin{bmatrix}
A&B\\
C&D
\end{bmatrix}\in M_2(\mathcal L(X))^H.$
\end{cor}

\begin{proof}
Clearly,
$M=\begin{bmatrix}
A&B\\
C&D
\end{bmatrix}=\begin{bmatrix}
0&1\\
1&0
\end{bmatrix}\begin{bmatrix}
D&C\\
B&A
\end{bmatrix}\begin{bmatrix}
0&1\\
1&0
\end{bmatrix},
$\\
where 
$$
M'=\begin{bmatrix}
D&C\\
B&A
\end{bmatrix}=\begin{bmatrix}
D&C\\
0&0
\end{bmatrix}+\begin{bmatrix}
0&0\\
B&A
\end{bmatrix}:=P+Q.
$$ 

We only need to prove $M'\in M_2(\mathcal L(X))^H.$

Compute that 
$$
P^H=\begin{bmatrix}
D^H&(D^H)^2C\\
0&0
\end{bmatrix}, Q^H=\begin{bmatrix}
0&0\\
(A^H)^2B&A^H
\end{bmatrix}.
$$

Obviously,
$$
P^HQ=\begin{bmatrix}
D^H&(D^H)^2C\\
0&0
\end{bmatrix}\begin{bmatrix}
0&0\\
B&A
\end{bmatrix}=0,
$$
$$
PQ^H=\begin{bmatrix}
D&C\\
0&0
\end{bmatrix}\begin{bmatrix}
0&0\\
(A^H)^2B&A^H
\end{bmatrix}=0,
$$
$$
Q^{\pi}PQP^{\pi}=\begin{bmatrix}
CBD^{\pi}&CA\\
0&0
\end{bmatrix}=0.
$$

According to Corollary 3.3, $M'$ has Hirano inverse, which follows that $M$ has Hirano inverse.
\end{proof}

\begin{exam}
Let $A=\begin{bmatrix}
1&1\\
0&0
\end{bmatrix}, B=\begin{bmatrix}
1&0\\
-1&0
\end{bmatrix}, C=\begin{bmatrix}
0&0\\
2&3
\end{bmatrix},\\D=\begin{bmatrix}
0&0\\
0&1
\end{bmatrix}$. Clearly, $A, D\in M_2(\mathbb C)^H.$ Observe that 
$$
AB=\begin{bmatrix}
1&1\\
0&0
\end{bmatrix}\begin{bmatrix}
1&0\\
-1&0
\end{bmatrix}=0,
$$
$$
BD^H=\begin{bmatrix}
1&0\\
-1&0
\end{bmatrix}\begin{bmatrix}
0&0\\
0&1
\end{bmatrix}=0,
$$
$$
D^{\pi}CB=\begin{bmatrix}
1&0\\
0&0
\end{bmatrix}\begin{bmatrix}
0&0\\
2&3
\end{bmatrix}\begin{bmatrix}
1&0\\
-1&0
\end{bmatrix}=0.
$$

By Theorem 3.7, $M=\begin{bmatrix}
A&B\\
C&D
\end{bmatrix}\in M_4(\mathbb C)^H.$
\end{exam}

\end{document}